\documentclass[12pt]{article}
\usepackage[a4paper]{geometry}
\usepackage{amssymb,amsbsy,graphicx}
\usepackage{times}
\newcommand{\cH}{{\cal H}}

\newcommand{\cQ}{{\cal Q}}
\newcommand{\cR}{{\cal R}}
\newcommand{\cS}{{\cal S}}
\newcommand{\cT}{{\cal T}}

\newcommand{\ZZ}{\mathbb{Z}}
\newcommand{\RR}{\mathbb{R}}
\newcommand{\NN}{\mathbb{N}}
\newcommand{\CC}{\mathbb{C}}

\newcommand{\Ab}{{\boldsymbol{A}}}
\newcommand{\Bb}{{\boldsymbol{B}}}
\newcommand{\Cb}{{\boldsymbol{C}}}

\newcommand{\Hb}{{\boldsymbol{H}}}

\newcommand{\Wb}{{\boldsymbol{W}}}

\newcommand{\ab}{{\boldsymbol{a}}}
\newcommand{\bb}{{\boldsymbol{b}}}
\newcommand{\cb}{{\boldsymbol{c}}}

\newcommand{\fb}{{\boldsymbol{f}}}

\newcommand{\hb}{{\boldsymbol{h}}}

\newcommand{\pb}{{\boldsymbol{p}}}

\newcommand{\vb}{{\boldsymbol{v}}}

\newcommand{\eop}{\hfill$\Box$}

\newtheorem{theorem}{Theorem}

\newtheorem{lemma}[theorem]{Lemma}
\newtheorem{corollary}[theorem]{Corollary}
\newtheorem{example}[theorem]{Example}
\newtheorem{remark}[theorem]{Remark}
\newtheorem{definition}[theorem]{Definition}

\newenvironment{pf}%
{\par\noindent\textbf{Proof:}~}%
{\eop\par\smallskip\par\noindent}

\newenvironment{pfof}[1]%
{\par\noindent\textbf{Proof of #1:}~}%
{\eop\par\smallskip\par\noindent}

\begin{document}
\title{Hermite subdivision schemes, exponential polynomial generation,
  and annihilators}
\author{Costanza Conti\thanks{DIEF, Universit\`a di Firenze, Viale
    Morgagni 40--44, I--50134 Firenze,
    Italy. \texttt{costanza.conti@unifi.it}} \and Mariantonia
  Cotronei\thanks{DIIES, Universit\`a Mediterranea di Reggio Calabria,
  Via Graziella, I--89122 Reggio Calabria,
  Italy. \texttt{mariantonia.cotronei@unirc.it}}
  \and Tomas Sauer\thanks{Lehrstuhl f\"ur Mathematik mit Schwerpunkt
    Digitale Bildverarbeitung, University of Passau, Innstr. 43,
    D--94032 Passau, Germany. \texttt{Tomas.Sauer@uni-passau.de}}}
\maketitle

\begin{abstract}
  We consider the question when the so--called \emph{spectral
    condition} for Hermite subdivision schemes extends to spaces
  generated by polynomials and exponential functions. The main tool
  are convolution operators that annihilate the space in question
  which apparently is a general concept in the study of various types
  of subdivision operators. Based on these annihilators, we
  characterize the spectral condition in terms of factorization of the
  subdivision operator.
\par\smallskip\noindent
{\bf Keywords} subdivision schemes; Hermite schemes; factorization; annihilators
 \par\smallskip\noindent
{\bf MSC} 65D15;  41A05; 42C15
	\end{abstract}

\section{Introduction}
Subdivision schemes are efficient iterative procedures based on the
repeated application of \emph{subdivision operators} which might
differ at different levels of iteration. Whenever convergent, they
generate functions that hopefully resemble the data used to start the
iterative procedure.

Subdivision operators act on bi-infinite sequences $\cb :
\ZZ \to \RR$ by means of a \emph{finitely} supported \emph{mask} $\ab :
\ZZ \to \RR$ in the convolution--like form
$$
\cS_\ab \cb = \sum_{\beta \in \ZZ} a( \cdot - 2\beta ) \, c(\beta).
$$
This type of operators has been generalized in various ways,
considering multivariate operators, operators with dilation factors
other than $2$  or
subdivision operators acting on vector or matrix data by means of matrix valued
masks. 
There is such a vast amount of literature
meanwhile that we do not even attempt to give references here.

It has been observed from very early on that preservation of
\emph{polynomial data} is an important property of subdivision
operators. For example, the preservation of constants, $\cS_\ab 1 =
1$, is necessary 
for the convergence of the subdivision schemes 
which iterate the same operator $\cS_\ab$.
More generally, the preservation of polynomial
spaces, $\cS_\ab \Pi_n = \Pi_n$, plays an important role in the
investigation of the
differentiability of the limit function of subdivision
schemes. In addition, there has been interest in  also preserving
functions other than polynomials, see for example 
\cite{Micchelli96}, and it is natural that such functions must be
exponential, i.e., of the form $e^{\lambda \cdot}$, cf. \cite{UnserBlu05}.

In
this paper we will consider preservation of such exponentials
by \emph{Hermite subdivision operators} which act on vector data but
with the particular understanding that these vectors represent function values
and consecutive derivatives up to a certain order. We will study the
preservation capability of such 
operators by means of a \emph{cancellation operator}, a concept that
applies to subdivision schemes in quite some generality. 
This is why, before we get to the main technical content of the paper,
we want to illustrate the idea and the concept through a few examples.

The simplest example deals with the preservation of constants, $\cS_\ab 1 = 1$.
Note that constant sequences are exactly the kernel of the
difference operator $\Delta$, defined as $\Delta \cb = c(\cdot + 1) -
c$; in other words: the difference operator is the simplest
\emph{cancellation operator} or \emph{annihilator} of the constant
functions. Now, whenever $\cS_\ab$ \emph{preserves} constants, then $\cS_\cb
= \Delta \cS_\ab$ 
is a subdivision operator that \emph{annihilates} the constants. As it can
easily be shown, any such operator can be written as $\cS_\cb = \cS_\bb
\Delta$ for some other finitely supported mask $\bb$, hence we get the
factorization $\Delta \cS_\ab = \cS_\bb \Delta$. Switching to the calculus of
\emph{symbols} which associates to a finitely supported sequence $\ab$
the \emph{Laurent polynomial}
$$
a^* (z) := \sum_{\alpha \in \ZZ} a(\alpha) \, z^\alpha,
$$
the factorization is equivalent to $(z^{-1}-1) a^* (z) = b^* (z) (z^{-2}-1)$
or, equivalently, to the famous ``zero at $\pi$'' condition $a^* (z) =
\left( z^{-1}+1 \right) \, b^* (z)$.

For a slightly more sophisticated example, suppose that now the
subdivision operator provides preservation of the subspace
\begin{equation}
  \label{eq:VdLambdaDef}
  V_{d,\Lambda} 
  = \mbox{span}\, \left\{
    1,x,\dots,x^p,e^{\lambda_1 x},e^{-\lambda_1 x},\dots,e^{\lambda_r
      x},e^{-\lambda_r x} \right\}, \qquad d = p+2r+1,
\end{equation}
in the sense that $\cS_\ab V^0_{d,\Lambda} \subseteq V^1_{d,\Lambda}$
where $V^j_{d,\Lambda}:=\left\{ v(2^{-j} \cdot) \;:\; v\in V_{d,\Lambda}\right\}$,
see, for example, \cite{ContiRomani11,UnserBlu05}. Again we
approach this problem in terms of cancellation, therefore determining
an operator $\cH_{d,\Lambda}$ such that $\cH_{d,\Lambda} V^0_{d,\Lambda} =
\{ 0 \}$.  Assuming that
$\cH_{d,\Lambda}$ is a \emph{convolution operator} (or \emph{LTI
  filter} in the language of signal processing, cf. \cite{Hamming98})
with \emph{impulse response} $\hb$, it is easily seen and well--known that
cancellation of the polynomials of degree at most $p$ implies that $\left( h^*
\right)^{(k)} (1) = 0$, $k=0,\dots,p$,
hence cancellation of the polynomial part of $V_{d,\Lambda}$ implies
that $h^* (z) = \left( z^{-1} - 1 \right)^{p+1} \, b_1^*
(z)$. Cancellation of an exponential sequence $e^{\lambda \cdot}$, on
the other hand, leads to
\begin{eqnarray*}
  0 & = & \sum_{j \in \ZZ} h(\cdot - j) e^{\lambda j} = \sum_{j \in
    \ZZ} h(j) e^{\lambda (\cdot - j)} = e^{\lambda \cdot} h^* (e^{-\lambda}),
\end{eqnarray*}
hence, the annihilation of the space implies that
$$
h^* (z) = b_2^* (z) \, \prod_{j=1}^r \left( z^{-1} - e^{\lambda_j}
\right) \left(z^{-1} - e^{-\lambda_j} \right).
$$
Summarizing, the simplest cancellation operator for $V_{d,\Lambda}$
takes the form
$$
h_{d,\Lambda}^* (z) = (z^{-1} - 1)^{p+1} \, \prod_{j=1}^r \left(
  z^{-1} - e^{\lambda_j} \right) \left(z^{-1} - e^{-\lambda_j} \right),
$$
and the associated factorization by means of cancellation operators
\begin{equation}
  \label{eq:HFactScalar}
  \cH_{d,2^{-1}\Lambda} \cS_\ab = \cS_\bb \cH_{d,\Lambda}
\end{equation}
is easily verified to be equivalent to the symbol factorization
\begin{equation}
  \label{eq:HFactScalarSymbol}
  a^* (z) = b^* (z) \, (z^{-1} + 1)^{p+1} \, \prod_{j=1}^r \left(
  z^{-1} + e^{\lambda_j/2} \right) \left(z^{-1} + e^{-\lambda_j/2} \right),
\end{equation}
given in \cite{UnserBlu05}.
Note that in (\ref{eq:HFactScalar}) one really has
to consider different spaces, hence a preservation property of the
form $\cS_\ab V^0_{d,\Lambda} \subseteq V^1_{d,\Lambda}$ because the result
of the subdivision operator corresponds to a sequence on the grid
$\ZZ/2$.

The last example considers \emph{Hermite subdivision schemes} which we
will investigate in more detail in the rest of this paper. In  Hermite
subdivision, the data
are vector valued sequences $\vb \in \ell^{d+1} (\ZZ)$ with the intuition
that the $k$--th component of such a sequence represents a $k$--th
derivative. Then, as considered for example in
\cite{ContiMerrienRomani14,DubucMerrien09,MerrienSauer12}, one
defines, for $f \in C^d (\RR)$, a sequence
$$
\vb_{f} \;:\; \alpha \mapsto \left[ f^{(j)}(\alpha) \;:\; j=0,\dots,d
\right], \qquad \alpha \in \ZZ,
$$
and asks when a subdivision operator $\cS_\Cb$ with \emph{matrix valued}
masks $\Cb \in \ell_{00}^{d \times d} (\ZZ)$ annihilates all $\vb_p$ for $p
\in \Pi_d$ which, by the aforementioned machinery, can again be used
to describe the \emph{spectral condition}, a ``polynomial
preservation'' rule introduced by Dubuc and Merrien in
\cite{DubucMerrien09}. Note  that it is no mistake or accident
that the letter $d$ appears for the maximal order of derivatives and
the maximal degree of polynomial cancellation -- the space dimension
and the order of derivatives are closely tied. It was then shown in
\cite{MerrienSauer12} that whenever $\cS_\Cb \vb_p = 0$ for all $p \in
\Pi_d$, then there exist a finitely supported $\Bb\in \ell_{00}^{(d+1)
  \times (d+1)} (\ZZ) $ such that
$$
C^* (z) = B^* (z) \, T_d^* (z^2), \qquad T_d^* (z) := \left[
  \begin{array}{cccc}
    z^{-1} - 1 & \frac12 & \dots & -\frac{1}{(d+1)!} \\
    & z^{-1} - 1 & \dots & -\frac{1}{d!} \\
    & & \ddots & \vdots \\
    & & & z^{-1} - 1
  \end{array}
\right].
$$
Since the operator $\cT$ acts for $f \in C^{d+1}$, $k=0,\cdots,d$, as
\begin{equation}
  \label{eq:TaylorFormula}
  \left( \cT v_f \right)_k (\alpha) = f^{(k)} (\alpha + 1) -
  \sum_{j=0}^{d-k} \frac{f^{(k+j)} (\alpha)}{j!} = f^{(d+1)} (\xi_k),
  \qquad \xi_k \in (\alpha,\alpha+1),
\end{equation}
hence measures the difference between a function and its Taylor
polynomial approximation at the neighboring point, it is
called the \emph{(complete) Taylor operator} of order $d$. That $\cT$
annihilates all $\vb_p$, $p \in \Pi_d$, is immediate from
(\ref{eq:TaylorFormula}).

It should have become clear by now that there is an obvious common structure
behind all these
examples. \emph{Preservation} of a subspace that can be
written as the kernel of a convolution operator is related to a
commuting property provided that the convolution operator factorize-s
or ``divides'' \emph{any} annihilator of the subspace. This can be
seen as a minimality property with respect to the partial ordering
given by divisibility and justifies the following terminology where we
identify any function $f \in V$ with the sequence $\fb= \left(
  f(\alpha) \;:\; \alpha \in \ZZ \right)$.

\begin{definition}
  A linear operator $\cH : \ell^m (\ZZ) \to \ell^m (\ZZ)$ is called a
  \emph{convolution operator} for a space $V$ if there exists a matrix
  sequence $\Hb\in \ell^{m\times m} (\ZZ)$, called the \emph{impulse
    response} of $\cH$, such that
  $$
  \cH \fb = \Hb*\fb = \sum_{\beta\in \ZZ} H (\cdot -\alpha) f(\alpha),
  \qquad \fb\in V.
  $$
\end{definition}

\begin{definition}
  A convolution operator $\cH : \ell^m (\ZZ) \to \ell^m (\ZZ)$ is called a
  \emph{minimal annihilator} for a space $V$ with respect to
 \begin{enumerate}
  \item \emph{subdivision}, if for any $\Cb \in \ell^{m \times m}_{00} (\ZZ)$
    such that $\cS_\Cb V = 0$ there exists $\Bb \in \ell^{m \times
      m}_{00} (\ZZ)$ with $\cS_\Cb = \cS_\Bb \cH$.
  \item \emph{convolution}, if for any $\Cb \in \ell^{m \times m}_{00} (\ZZ)$
    such that $\Cb * V = \{0\}$ there exists $\Bb \in \ell^{m \times
      m}_{00} (\ZZ)$ with $\Cb = \Bb * \Hb$,
  \end{enumerate}
  respectively. If $\cH$ satisfies both properties it is simply called
  a \emph{minimal annihilator}.
\end{definition}

The goal of this paper is to use this general concept to understand
preservation of exponentials and polynomials by  Hermite subdivision
schemes where the subdivision operators will have to vary with the
iteration level; some call this \emph{nonstationary}, some
\emph{nonuniform} operators, but the problem is too 
interesting to dwell on such niceties here and therefore we
omit it as the name of a property that is not fulfilled anyway is simply
irrelevant. 

In more technical terms, we will derive the analogy of the Taylor operator for
the case of preservation of exponentials and prove in
Theorem~\ref{thm:ConvFactorization} that it is again a minimal
annihilator. We will see that even the cancellation operator depends
only on the space $V_{d,\Lambda}$ and on the level.
We will also see that the existence of the annihilator operator is strongly connected with the factorization  of the subdivision operator satisfying specific preservation properties. 

The organization of the paper is as follows. 
After providing the necessary notation and terminology, the main results on Hermite subdivision schemes and their reproduction capabilities
will be derived in Section~\ref{sec:hermite}.
To better explain the underlying ideas, we will first consider the case of adding a single frequency to the polynomial space and then extend the results and methods to an arbitrary number of frequencies. These descriptions will be in terms of appropriate cancellation operators.
Thereafter, in Section~\ref{sec:factorization} we will use these cancellation operators to derive factorization properties which will also verify that the cancellation operators are minimal.
Finally, we will illustrate our  results with specific examples.

\section{Subdivision schemes and notation}
We begin by fixing the notation and 
recalling some known facts about subdivision schemes. 
We denote by $\ell^m \left( \ZZ \right)$ and $\ell^{m \times m} \left(
  \ZZ \right)$ the linear spaces of all sequences of  $m$--vectors and
$m \times m$ matrices, respectively. 
Operators acting on that spaces are denoted by capital calligrafic
letter. Sequences in $\ell^m \left( \ZZ \right)$ and $\ell^{m \times
  m} \left( \ZZ \right)$ will be denoted by boldface lower case and
upper case letters, respectively. 
In particular, $\cb\in \ell^m \left(\ZZ \right)$ is $\cb= \left( c(\alpha)
  \;:\; \alpha \in \ZZ \right)$, while $\Ab\in \ell^{m \times m}
\left(\ZZ \right)$ stands for $\Ab = \left( A(\alpha) \;:\, \alpha \in
  \ZZ \right)$, indexing $A\in \RR^{m\times m}$ as 
$A=\left[
  a_{jk}:\, 	j,k=0,\dots,m-1\right]$.
As usual,  $\ell_{00}^m \left( \ZZ \right)$ and
$\ell_{00}^{m\times m} \left( \ZZ \right)$ will denote the subspaces
of finitely supported  sequences, and $\NN_0$ denotes the set $\{0,1,2,\ldots\}$.

For $\Ab\in \ell_{00}^{m\times m} \left( \ZZ \right)$ and
$\cb\in \ell_{00}^m \left( \ZZ \right)$  we define the associated
\emph{symbols} as the \emph{Laurent polynomials}
$$
 A^* (z):=\sum_{\alpha \in \ZZ} A(\alpha) \, z^{\alpha}, \qquad
 c^* (z):=\sum_{\alpha \in \ZZ} c(\alpha) \,
 z^{\alpha}, \qquad z \in \CC \setminus \{0\}.
$$
For $\Ab\in \ell_{00}^{m\times r} \left( \ZZ \right)$ and $\Bb\in
\ell_{00}^{r\times q} \left( \ZZ \right)$ the \emph{convolution}
$\Cb=\Ab*\Bb$ in  $\ell_{00}^{m\times q} \left( \ZZ \right)$ is defined
as usually as
$$
 C(\alpha):=\sum_{\beta \in \ZZ} A(\beta)B(\alpha-\beta), \qquad \alpha\in \ZZ.
$$
The \emph{subdivision operator} $\cS_\Ab : \ell^m(\ZZ) \rightarrow
\ell^m(\ZZ)$ based on the matrix sequence or
\emph{mask} $\Ab \in \ell^{m \times m}_{00}(\ZZ)$ is defined as
\begin{equation}\label{def:sub_operator}
  \cS_\Ab \cb (\alpha) = \sum_{\beta \in \ZZ} A \left(\alpha -2
    \beta \right) c(\beta), \qquad \alpha \in \ZZ,\qquad  \hbox{for
    all}\ \cb \in \ell^m \left( \ZZ \right).
\end{equation}
Alternatively, using  symbol calculus notation, we can also decribe the
action of the subdivision operator in the form 
\begin{equation}
  \label{lab1} \left( \cS_\Ab \cb \right)^* (z) = A^* (z) \, c^*
  \left( z^2 \right), \qquad z \in  \CC \setminus \{0\},
\end{equation}
though, in strict formalism, (\ref{lab1}) is only valid for  $\cb \in \ell_{00}^m (\ZZ)$.

A \emph{subdivision scheme} consists of the successive application of
potentially different subdivision operators $\cS_{\Ab^{[n]}}$,
constructed from a sequence of masks $\left( \Ab^{[n]} \;:\; n \in
  \NN_0 \right)$ where $\Ab^{[n]}=\left( A^{[n]} (\alpha) \;:\; \alpha
  \in \ZZ \right) \in \ell_{00}^{m \times m}(\ZZ)$ is called the
\emph{level $n$ subdivision mask} and is assumed to be of finite
support.  
Accordingly, a sequence of matrix \emph{symbols}
$\left( (A^{[n]})^*(z) \;:\; n \in \NN_0 \right)$ characterizes such schemes.

For some initial sequence $\cb^{[0]} \in \ell^m (\ZZ)$ the subdivision
scheme iteratively produces sequences
$$
\cb^{[n+1]} := \cS_{\Ab^{[n]}}\cb^{[n]}, \qquad n \in \NN_0,
$$
whose elements can be interpreted as function values at
$2^{-n-1}\ZZ$, from which one can define convergence in the usual way.

\section{Hermite subdivision schemes and reproduction}
\label{sec:hermite}
As already mentioned, Hermite subdivision schemes act on vector valued
data $\cb \in \ell^{d+1} (\ZZ)$, whose $k$-th component 
corresponds to a $k$--th derivative.
We are interested in studying the exponential and polynomial
preservation capabilities of such kind of
schemes. 

A preliminary simple observation is that for $f\in C^d(\RR)$ and for
$g:=f(2^{-n} \cdot)$ we clearly have
$\frac{d^r}{dx^r} g = 2^{-nr} \frac{d^r f}{dx^r}f (2^{-n}\cdot)$, $r=0,\cdots,d$.
Hence
\begin{equation}\label{eq:1}
 \left[ \frac {d^j}{dx^j} f (2^n \cdot):\, j=0,\dots,d\right]= D^n
\left[ f^{(j)} (2^{-n}\cdots):\, j=0,\dots,d\right] ,
\end{equation}
where
$$ 
D = \left[
  \begin{array}{ccccc}
    1& 0 & 0 &\cdots & 0\\ 0& \frac12 &0&\cdots & 0\\ \vdots \\ 0& 0&
    0& \cdots& \frac{1}{2^d}
  \end{array}\right]. 
$$
Since the sequence $\cb^{[n]}$ is related to evaluations on the grid $2^{-n}\ZZ$,
we consider Hermite subdivision schemes with the $n$-th iteration of
the following type: 
\begin{equation}
  \label{eq:HermSubd0}
  D^{n+1} \cb^{[n+1]}=\sum_{\beta\in \ZZ} A^{[n]}(\cdot -2\beta) D^n c^{[n]}(\beta),
\end{equation}
where in ``usual'' Hermite subdivision the mask is the same over all
levels, i.e., $\Ab^{[n]} = \Ab$, $n \in \NN_0$.
Setting $\widetilde \Ab^{[n]} := D^{-n-1} \Ab^{[n]} D^n$,
(\ref{eq:HermSubd0}) fits into the framework of Section~2 with the $n$-th
subdivision operator of the form
\begin{equation}\label{eq:HermSubd}
\cb^{[n+1]} = S_{\widetilde A^{[n]}} \cb^{[n]}
= \sum_{\beta\in \ZZ} \widetilde A^{[n]}(\cdot -2\beta) \, c^{[n]}(\beta).
\end{equation}

\subsection{Single exponential frequency}
\label{subsec:singlefreq}
In the first step of our analysis of the stepwise reproduction
capability of a Hermite subdivision scheme of type
(\ref{eq:HermSubd0}), we add only \emph{a single} pair of exponentials
$e^{\pm \lambda x}$ and consider the space
\begin{equation}
  \label{eq:VlambdaDef}
  V_{d,\lambda}  = \mbox{\rm span}
  \left\{1,x,\ldots,x^{d-2},e^{\lambda x},e^{-\lambda x}\right\},
  \qquad \lambda\in \RR \cup i\RR. 
\end{equation}
To keep notation simple and to better explain the underlying ideas, we will
first carefully investigate this situation and then extend it in a
quite straightforward fashion to the general case.

\begin{remark}
  As can be seen in (\ref{eq:VlambdaDef}), the addition of an
  exponential frequency $\lambda$ always means the addition of the
  pair $e^{\pm \lambda \cdot}$ of functions. On the one hand, this is
  motivated by the fact that choosing $\lambda = i$ equals
  reproduction of the trigonometric functions $\sin x$ and $\cos
  x$. Moreover, our approach to construct the annihilator and the
  factorization actually strongly depends on the presence of this pair
  of frequencies. Whether or not similar results will be available for
  the case where only $e^{\lambda \cdot}$ but not $e^{-\lambda
    \cdot}$, we do not know at present.
\end{remark}

\noindent
For any function $f\in C^d(\RR)$ and any integer $n \in \NN_0$ we
consider the two vector sequences ${\widetilde \vb}_{f,{n}}$,
$\vb_{f,{n}}\in \ell^{d+1}(\ZZ)$, defined, for $\alpha \in \ZZ$, as
$$
\widetilde v_{f,{n}} (\alpha) := \left[
  \begin{array}{c}
    f( 2^{-n} \alpha ) \\ f'( 2^{-n} \alpha )\\ \vdots \\ f^{(d)}(
    2^{-n} \alpha )
  \end{array}
\right],
\qquad v_{f,{n}} (\alpha) :=D^n \widetilde v_{f,{n}} (\alpha) =
\left[
  \begin{array}{c}
    f( 2^{-n} \alpha ) \\ 2^{-n} f'( 2^{-n} \alpha )\\ \vdots \\
    2^{-nd} f^{(d)}(2^{-n} \alpha )
  \end{array}
\right].
$$
We simply write $\vb_{f}=\widetilde \vb_{f}$ when $n=0$.

\begin{definition}\label{def:spectrVd}
  A mask $\Ab^{[n]}\in \ell_{00}^{(d+1)\times(d+1)}(\ZZ)$ or its
  associated subdivision operator
  $\cS_{{\Ab}^{[n]}}$ satisfies the $V_{d,\lambda}$-spectral condition
  if 
  $$
  \cS_{{\Ab}^{[n]}} \vb_{f,n} = \vb_{f,{n+1}},\qquad  f \in
  V_{d,\lambda}.
  $$
  Equivalently, the mask $\widetilde \Ab^{[n]} = D^{-(n+1)} \Ab^{[n]}
  D^n$ satisfies the $V_{d,\lambda}$-spectral condition if
  $$
  \cS_{\widetilde \Ab^{[n]}} {\widetilde \vb}_{f,n} = {\widetilde
    \vb}_{f,{n+1}}, \qquad f\in V_{d,\lambda}.
  $$
\end{definition}

\begin{remark}
  It is important to observe that Definition \ref{def:spectrVd} is
  fully consistent with \cite[Definition 1]{MerrienSauer12} though
  formulated in a slightly different way taking into account the
  stronger form of level dependency needed for the reproduction of
  exponentials.
\end{remark}

\noindent
Since we plan to extend difference operators and Taylor operators, we
next recall their formal definitions.


\begin{definition}\label{def:TaylorOp}
  The Taylor operator $\cT_d$ of order $d$, acting on
  $\ell^{d+1}(\ZZ)$ is defined as 
  \begin{equation}\label{eq:TaylorOp}
    \cT_d:=\left[ \begin{array}{ccccc}
        \Delta&-1 & \cdots &-\frac 1{(d-1)!} &-\frac 1{d!}\\
        & \Delta &\ddots & \vdots & \vdots\\
        &   &\ddots & -1& \vdots\\
        &   &  & \Delta& -1\\
        &   &  & & \Delta
      \end{array}
    \right],
  \end{equation}
	where $\Delta$ is the forward difference operator.
\end{definition}

The symbol of the Taylor operator then takes the form
\begin{equation}\label{eq:SymbTaylorOp}
  T_d^*(z):=\left[ \begin{array}{ccccc}
      (z^{-1}-1) &-1 & \cdots &-\frac 1{(d-1)!} &-\frac 1{d!}\\
      & (z^{-1} - 1) &\ddots & \vdots & \vdots\\
      &   &\ddots & -1& \vdots\\
      &   &  & (z^{-1} - 1)& -1\\
      &   &  & & (z^{-1}-1)
    \end{array}
  \right].
\end{equation}

\begin{definition}\label{def:AnnOp}
  A \emph{level-$n$ cancellation  operator} $\cH^{[n]} :
  \ell^{d+1}(\ZZ) \to \ell^{d+1}(\ZZ)$ for a linear function space $V
  \subset C^d (\RR)$ 
  is a convolution operator such that
  \begin{equation}\label{eq:AnnOp}
    \cH^{[n]} \vb_{f,n} = \sum_{\alpha\in \ZZ} H^{[n]} (\cdot-\alpha)
    v_{f,n}(\alpha)=0, \qquad f\in V. 
  \end{equation}
  By $\cH^{[n]}_{d,\lambda}$ we denote a level-$n$ cancellation
  operator for the function space spanned by $V_{d,\lambda}$.
\end{definition}

\begin{lemma}\label{lem:AnnihilOp}
  An operator $\cH^{[n]}_{d,\lambda}$ is a level-$n$ cancellation
  operator for the space $V_{d,\lambda}$ if it satisfies
  \begin{equation}
    \label{eq:AnnihilOp1}
    \left( H^{[n]}_{d,\lambda} \right)^*(z) = \left[
      \begin{array}{cc}
        T^*_{d-2}(z) & * \\
        0 & *
      \end{array}
    \right]
  \end{equation}
  and
  \begin{equation}
    \label{eq:AnnihilOp2}
    \left( H^{[n]}_{d,\lambda} \right)^* \left( e^{\mp 2^{-n} \lambda} \right) D^n
    \left[
      \begin{array}{c}
        1 \\ \pm \lambda \\ \vdots \\ (\pm \lambda)^d
      \end{array}
    \right]
    = 0.
    \end{equation}
\end{lemma}

\begin{pf}
  To annihilate polynomials of degree $d-2$,
  condition (\ref{eq:AnnOp}) has to be satisfied for the vector sequences
  $$
  \vb_{(\cdot)^j,n} = D^n\left[ (2^{-n} \cdot )^j, j (2^{-n}
    \cdot)^{j-1},\dots, j!,
    \underbrace{0,\,0, \cdots,\,0}_{d-j}\right]^T,\qquad
  j=0,\ldots,d-2,
  $$
  and since these sequences are exactly annihilated by the complete
  Taylor operator as shown in \cite{MerrienSauer12}, any decomposition
  of the form (\ref{eq:AnnihilOp1}) annihilates polynomials of degree
  at most $d-2$.

  To describe cancellation of exponentials, we first observe that
  $$
  \vb_{e^{\pm \lambda \cdot},n} = e^{\pm \lambda 2^{-n} \cdot} D^n
  \left[
    \begin{array}{c}
      1 \\ \pm \lambda \\ \vdots \\ (\pm \lambda)^d
    \end{array}
  \right] ,
  $$
  so that the condition becomes
  \begin{eqnarray*}
      0 & = & \sum_{\alpha \in \ZZ} H_{d,\lambda}^{[n]} ( \cdot - \alpha
      ) \, e^{\pm 2^{-n} \lambda \alpha} D^n \left[
        \begin{array}{c}
          1 \\ \pm \lambda \\ \vdots \\ (\pm \lambda)^d
        \end{array}
    \right]  \\
    & = & \sum_{\alpha \in \ZZ} H_{d,\lambda}^{[n]} ( \alpha ) \,
    e^{\pm 2^{-n} \lambda (\cdot - \alpha)} D^n\left[
      \begin{array}{c}
        1 \\ \pm \lambda \\ \vdots \\ (\pm \lambda)^d
      \end{array}
    \right] \\
    & = & e^{\pm 2^{-n} \lambda \cdot} \left( H_{d,\lambda}^{[n]}
    \right)^* ( e^{\mp 2^{-n} \lambda} ) \, D^n \left[
      \begin{array}{c}
        1 \\ \pm \lambda \\ \vdots \\ (\pm \lambda)^d
      \end{array}
    \right]
  \end{eqnarray*}
  which yields (\ref{eq:AnnihilOp2}).
\end{pf}

\begin{remark}\label{rem_marianto}
  If we are able to find an operator $\cH_{d,\lambda}$ that satisfies
  (\ref{eq:AnnihilOp1}) and
  \begin{equation}\label{eq:AnnihilOp2n0}
    H_{d,\lambda}^* \left( e^{\mp \lambda} \right)
    \left[
      \begin{array}{c}
        1 \\ \pm \lambda \\ \vdots \\ (\pm \lambda)^d
      \end{array}
    \right]
    = 0,
  \end{equation}
  then we automatically obtain level-$n$ cancellation operators
  $\cH^{[n]}_{d,\lambda}$ for any $n \in \NN_0$ by setting
  $$
  \cH^{[n]}_{d,\lambda}=\cH_{d,2^{-n}\lambda}.
  $$
  In fact, this follows from the simple observation that the identity
  $$
  H_{d,2^{-n}\lambda}^* \left( e^{\mp 2^{-n} \lambda} \right) D^n
  \left[
    \begin{array}{c}
      1 \\ \pm \lambda \\ \vdots \\ (\pm \lambda)^d
    \end{array}
  \right]
  = H_{d,2^{-n}\lambda}^*\left( e^{\mp 2^{-n} \lambda} \right)
  \left[
    \begin{array}{c}
      1 \\ \pm \frac{\lambda}{2^n} \\ \vdots \\ (\pm \frac{\lambda}{2^n})^d
    \end{array}
  \right]
  = 0,
  $$
  is equivalent to (\ref{eq:AnnihilOp2n0}), as can be verified by
  just replacing $\lambda$ with $2^{-n} \lambda$.
\end{remark}
				
\noindent 		
In view of Remark~\ref{rem_marianto} we see that to generate a
level-$n$ cancellation operator we just need to construct a level-$0$
cancellation operator. Therefore we continue with the analysis of
$\cH^{[0]}_{d,\lambda}$ which will be simply denoted by $\cH_{d,\lambda}$.

The next step is now to construct a cancellation operator which will
eventually even turn out to be a minimal one.
\smallskip  \noindent  Based on Lemma~\ref{lem:AnnihilOp}, the structure of the
cancellation operator $\cH_{d,\lambda}$ for the space $V_{d,\lambda}$
can now be derived. Indeed, we write its symbol in the general form
\begin{equation}\label{eq:Hcase1}
H_{d,\lambda}^*(z) = \left[
  \begin{array}{cc}
    T_{d-2}^*(z) & Q^*(z) \\
    0 & R^*(z)
  \end{array}
\right],\quad Q^*(z)\in \RR^{(d-2)\times 2},\, R^*(z)\in \RR^{2\times 2}
\end{equation}
and determine the remaining part of $H_{d,\lambda}^*(z)$, namely the Laurent
polynomial matrices $Q^*(z)$ and $R^*(z)$. 
To this aim, we begin to explicitly compute  the first line
$(H^*_{d,\lambda})_{0,:}(z)$, where the ``$:$'' is to be understood in
the sense of Matlab notation.

\begin{lemma}\label{lem:Line1Compl}
  The condition
  \begin{equation}
    \label{eq:Line1ComplCond}
    (H^*_{d,\lambda})_{0,:} \left( e^{\mp \lambda} \right)
		\left[
      \begin{array}{c}
        1 \\ \pm \lambda \\ \vdots \\ (\pm \lambda)^d
      \end{array}
    \right] = 0\,,
  \end{equation}
  can be fulfilled by setting
  \begin{equation}
    \label{eq:Line1ComplCond0d-1}
    (H_{d,\lambda})_{0,d-1} = h_{0,d-1} = \frac{\lambda^{1-d}}2 \left\{
      \begin{array}{lcl}
        e^{-\lambda} - e^\lambda + 2\displaystyle\sum_{2j+1 \le d-2}
        \frac{\lambda^{2j+1}}{(2j+1)!}, & \; & d \in 2 \ZZ, \\
        -\left(e^{-\lambda} + e^\lambda - 2\displaystyle\sum_{2j \le d-2}
        \frac{\lambda^{2j}}{(2j)!}\right), & \; & d \in 2 \ZZ+1,
      \end{array}
    \right.
  \end{equation}
  and
  \begin{equation}
    \label{eq:Line1ComplCond0d}
    (H_{d,\lambda})_{0,d} = h_{0,d} = \frac{\lambda^{-d}}2 \left\{
      \begin{array}{lcl}
        -\left(e^{-\lambda} + e^\lambda + 2\displaystyle\sum_{2j \le d-2}
        \frac{\lambda^{2j}}{(2j)!}\right),  & \; & d \in 2 \ZZ, \\
        e^{-\lambda} - e^\lambda + 2\displaystyle\sum_{2j+1 \le d-2}
        \frac{\lambda^{2j+1}}{(2j+1)!},& \; & d \in 2 \ZZ+1,
      \end{array}
    \right.
  \end{equation}
\end{lemma}

\begin{pf}
  Due to (\ref{eq:SymbTaylorOp}),
  the identity (\ref{eq:Line1ComplCond}) can be written as
  \begin{eqnarray*}
    0 & = & e^{\pm \lambda} - 1 - \sum_{j=1}^{d-2} \frac{(\pm
      \lambda)^k}{k!} + (\pm \lambda)^{d-1} h_{0,d-1} + (\pm
    \lambda)^d h_{0,d} \\
    & = & e^{\pm \lambda} - t_{d-2} \left[ e^{\pm \lambda \cdot}
    \right] (1) +  (\pm \lambda)^{d-1} h_{0,d-1} + (\pm
    \lambda)^d h_{0,d},
  \end{eqnarray*}
  where
  $$
  t_k [f] = \sum_{j=0}^k \frac{f^{(j)}(0)}{j!} (\cdot)^j\,,
  $$
  denotes the Taylor polynomial of $f$ of order $k$
  expanded at
  $0$. Adding and subtracting the above conditions we get
  \begin{eqnarray*}
    0 & = & \left( e^{\lambda} \pm e^{-\lambda} \right) - t_{d-2} \left[
      e^{\lambda \cdot} \pm e^{-\lambda \cdot} \right] (1) \\
    & & \quad + \left(
    \lambda^{d-1} \pm (-\lambda)^{d-1} \right) h_{0,d-1} + \left(
    \lambda^{d} \pm (-\lambda)^{d} \right) h_{0,d}.
  \end{eqnarray*}
  If $d$ is even, this implies that
  \begin{eqnarray*}
    h_{0,d-1} & = & \frac{e^{-\lambda} - e^\lambda - t_{d-2} \left[
        e^{-\lambda \cdot} - e^{\lambda \cdot} \right] (1)}{2
      \lambda^{d-1}}, \\
    h_{0,d} & = & -\frac{e^{-\lambda} + e^\lambda - t_{d-2} \left[
        e^{-\lambda \cdot} + e^{\lambda \cdot} \right] (1)}{2
      \lambda^d},
  \end{eqnarray*}
  while for odd $d$ we get
  \begin{eqnarray*}
    h_{0,d-1} & = &- \frac{e^{-\lambda} + e^\lambda - t_{d-2} \left[
        e^{-\lambda \cdot} + e^{\lambda \cdot} \right] (1)}{2
      \lambda^{d-1}}, \\
    h_{0,d} & = & \frac{e^{-\lambda} - e^\lambda - t_{d-2} \left[
        e^{-\lambda \cdot} - e^{\lambda \cdot} \right] (1)}{2
      \lambda^d}.
  \end{eqnarray*}
  Since
  $$
  \frac{d^k}{dx^k} \left( e^{\lambda x} \pm e^{-\lambda x} \right) =
  \lambda^k \left( e^{\lambda x} \pm (-1)^k e^{-\lambda x} \right),
  $$
  we have that
  $$
  t_{d-2} \left[ e^{\lambda \cdot} - e^{-\lambda \cdot} \right]
  (1)
  = 2\lambda + \frac23 \lambda^3 + \cdots = 2 \sum_{2j+1 \le d-2}
  \frac{\lambda^{2j+1}}{(2j+1)!}\,,
  $$
  and
  $$
  t_{d-2} \left[ e^{\lambda \cdot} + e^{-\lambda \cdot} \right]
  (1) = 2 + \lambda^2 + \cdots = 2 \sum_{2j \le d-2}
  \frac{\lambda^{2j}}{(2j)!}.
  $$
Substituting these identities readily gives (\ref{eq:Line1ComplCond0d-1}) and (\ref{eq:Line1ComplCond0d}).
\end{pf}

\smallskip \noindent Taking into account the structure of $H_{d,\lambda}^*(z)$, we can now easily
give also the entries of the other lines.\\

\begin{corollary}
  For $k=0,\dots,d-2$, we have that
  \begin{eqnarray}
    \label{eq:hkd-1Formula}
    h_{k,d-1} & = & \frac{\lambda^{1-d+k}}2 \left\{
      \begin{array}{lcl}
        e^{-\lambda} - e^\lambda + 2\displaystyle\sum_{2j+1 \le d-2-k}
        \frac{\lambda^{2j+1}}{(2j+1)!}, & \quad & d-k \in 2 \ZZ, \\
       -\left(e^{-\lambda} + e^\lambda - 2\displaystyle\sum_{2j \le d-2-k}
        \frac{\lambda^{2j}}{(2j)!}\right), & \quad & d-k \in 2 \ZZ+1,
      \end{array}
    \right. \\
    \label{eq:hkdFormula}
    h_{k,d} & = & \frac{\lambda^{-d+k}}2 \left\{
    \begin{array}{lcl}
      -\left(e^{-\lambda} + e^\lambda - 2\displaystyle\sum_{2j \le d-2-k}
      \frac{\lambda^{2j}}{(2j)!}\right),  & \quad & d-k \in 2 \ZZ, \\
      e^{-\lambda} - e^\lambda + 2\displaystyle\sum_{2j+1 \le d-2-k}
      \frac{\lambda^{2j+1}}{(2j+1)!},& \quad & d-k \in 2 \ZZ+1,
    \end{array}
  \right.
  \end{eqnarray}
  in particular, $h_{k-1,d-1} = h_{k,d}$.
\end{corollary}

\noindent
To complete the construction of $H_{d,\lambda}^*(z)$, we have to define the
lower right block $R^*(z)$ as
\begin{eqnarray}
  \label{eq:RzDefinition}
  R^* (z)& = &\left[
    \begin{array}{cc}
      z^{-1} - \displaystyle\frac{e^\lambda + e^{-\lambda}}{2} &
      \displaystyle\frac{e^{-\lambda} - e^\lambda}{2\lambda} \\[5mm]
      \lambda \displaystyle\frac{e^{-\lambda} - e^\lambda}{2} & z^{-1} -
      \displaystyle\frac{e^\lambda + e^{-\lambda}}{2}
    \end{array}
  \right]\\
	&=&L_{d,\lambda}\,\left[\begin{array}{cc}
	z^{-1}-e^{\lambda} & 0\\
	0 & z^{-1}-e^{-\lambda} \end{array}\right] \, L_{d,\lambda}^{-1},
\end{eqnarray}
where
$$L_{d,\lambda}=\left[\begin{array}{cc}
	\lambda^{d-1} & (-\lambda)^{d-1} \\
	 \lambda^{d} & (-\lambda)^{d} \end{array}\right],
	$$
for which the validity of (\ref{eq:AnnihilOp2}) is easily verified by
direct computations.

\begin{example}
As an example, we  provide the explicit structures of $\cH_{2,\lambda}$, $\cH_{3,\lambda}$ for the spaces $V_{2,\lambda}=\mbox{\rm span}\left\{1,e^{-\lambda x},e^{\lambda x}\right\}$ and
$V_{3,\lambda}=\mbox{\rm span}\left\{1,x,e^{-\lambda x},e^{\lambda x}\right\}$:

\begin{equation}\label{esempio:H_2}
H_{2,\lambda}^*(z)=\left[\begin{array}{ccc}
z^{-1}-1 & \displaystyle{\frac {e^{-\lambda}-e^\lambda}{2\lambda}}& -\displaystyle{\frac {e^{-\lambda}+e^\lambda-2}{2\lambda^2}} \\[5mm]
0&\displaystyle{z^{-1}-\frac {e^{-\lambda}+e^{\lambda}}2}&  \displaystyle{\frac {e^{-\lambda}-e^{\lambda}}{2\lambda}   } \\[5mm]
0&\displaystyle{\lambda\frac {e^{-\lambda}-e^{\lambda}}{2}   } & \displaystyle{z^{-1}-\frac {e^{-\lambda}+e^{\lambda}}2}
\end{array}
\right]\,,
\end{equation}

and

\begin{equation}\label{esempio:H_3}
H_{3,\lambda}^*(z)=\left[\begin{array}{cccc}
z^{-1}-1 & -1 & \displaystyle{\frac{2-e^{-\lambda}-e^{\lambda}}{2\lambda^2}} & \displaystyle{\frac{2\lambda+e^{-\lambda}-e^{\lambda}}{2\lambda^3}}\\[5mm]
0 & z^{-1}-1 & \displaystyle{\frac {e^{-\lambda}-e^\lambda}{2\lambda}}& -\displaystyle{\frac {e^{-\lambda}+e^\lambda-2}{2\lambda^2}} \\[5mm]
0&0&\displaystyle{z^{-1}-\frac {e^{-\lambda}+e^{\lambda}}2}&  \displaystyle{\frac {e^{-\lambda}-e^{\lambda}}{2\lambda}   } \\[5mm]
0&0&\displaystyle{\lambda\frac {e^{-\lambda}-e^{\lambda}}{2}   } & \displaystyle{z^{-1}-\frac {e^{-\lambda}+e^{\lambda}}2}
\end{array}
\right]$$
$$=\left[\begin{array}{cccc}
z^{-1}-1 & -1 & \displaystyle{\frac{2-e^{-\lambda}-e^{\lambda}}{2\lambda^2}} & \displaystyle{\frac{2\lambda+e^{-\lambda}-e^{\lambda}}{2\lambda^3}}\\[5mm]
0 & & & \\[5mm]
0&&H^*_{2,\lambda}(z)&   \\[5mm]
0&& &
\end{array}
\right]\,.
\end{equation}
\end{example}

\noindent
Of course, the above construction of $\cH_{d,\lambda}$ is only one of
many possibilities to construct a cancellation operator for
$V_{d,\lambda}$. However, our construction is well--chosen in the sense that it
includes the Taylor operator as action on the polynomials and that it
in fact naturally extends the Taylor operator. 

\begin{theorem}\label{thm:HdlconvTd}
  \begin{equation}
    \label{eq:HdlconvTd}
    \lim_{\lambda \to 0} \cH_{d,\lambda} = \cT_d.
  \end{equation}
\end{theorem}

\begin{pf}
  It is obvious from (\ref{eq:RzDefinition}) that
  $$
  R^* (z) \to \left[
    \begin{array}{cc}
      z^{-1} - 1 & -1 \\ 0 & z^{-1} - 1
    \end{array}
  \right]\,,
  $$
  as $\lambda \to 0$, hence it suffices to show that $h_{k,d-1} \to
  -\frac{1}{(d-1-k)!}$ and $h_{k,d} \to -\frac{1}{(d-k)!}$ as
  $\lambda \to 0$. Suppose that $d-k$ is even in which case we get
  \begin{eqnarray*}
    h_{k,d-1} & = & \frac{e^{-\lambda} - e^\lambda - t_{d-2-k} \left[
        e^{-\lambda \cdot} - e^{\lambda \cdot} \right](1)}{2
      \lambda^{d-1-k}}
    = \frac{1}{2 \lambda^{d-1-k}} \sum_{j=d-1-k}^\infty  \left( (-1)^j
      - 1 \right) \frac{\lambda^j}{j!} \\
    & = & -\frac{1}{(d-1-k)!} +
    \lambda^2 \sum_{j=d+1-k} \frac{(-1)^j + 1}{2 j!} \lambda^{j-(d+1-k)},
  \end{eqnarray*}
  which converges as desired when $\lambda \to 0$. The arguments for
  $h_{k,d}$ and the case of odd $d-k$ are identical.
\end{pf}

\subsection{Multiple exponential frequencies}
\label{subsec:multiplefreq}
Having understood the case of a single frequency $\lambda$, it is not
hard any more to extend the construction to arbitrary sets of frequencies. To
that end, let $\Lambda = \left\{ \lambda_1,\dots,\lambda_r
\right\}$ consist of $r$ different frequencies, all either real or
purely imaginary, and let us construct a cancellation operator
$\cH_{d,\Lambda}$ for the space
$$
V_{d,\Lambda} := \mbox{span\,} \left\{ 1,\dots,x^p,e^{\pm \lambda_1
    \cdot}, \dots, e^{\pm \lambda_r \cdot} \right\}, \qquad d = p+2r.
$$
The conditions for cancellation extend in a straightforward way.

\begin{lemma}\label{lem:HLambdaCancel}
  The operator $\cH_{d,\Lambda}$ with symbol
  $$
  H_{d,\Lambda}^*(z) = \left[
    \begin{array}{cc}
      T^*_p(z) &  Q^*(z) \\ 0 &  R^*(z)
    \end{array}
  \right],\quad Q^*(z)\in \RR^{p\times 2r},\, R^*(z)\in \RR^{2r\times 2r},
  $$
  annihilates $V_{d,\Lambda}$ if and only if
  \begin{equation}
    \label{eq:HLambdaCancel}
    H_{d,\Lambda}^*(z) \left( e^{\mp \lambda_j} \right) \left[
      \begin{array}{c}
        1 \\ \pm \lambda_j \\ \vdots \\ (\pm \lambda_j)^d
      \end{array}
    \right] = 0, \qquad j=1,\dots,r.
  \end{equation}
\end{lemma}

\begin{pf}
  Since the Taylor part of $\cH_{d,\Lambda}$ annihilates the   polynomials, we only need to perform the computations used to derive (\ref{eq:AnnihilOp2}) for any $\lambda_j$ to show that cancellation of the exponential polynomials is equivalent to
  (\ref{eq:HLambdaCancel}).
\end{pf}

\medskip \noindent
The construction of $\cH_{d,\Lambda}$ now follows the same lines as
 before, namely by determining the matrix symbol $Q^*(z)$.
 For the first row we now get, for $j=1,\dots,r$, the conditions
\begin{eqnarray*}
  0 & = & e^{\pm \lambda_j} - 1 - \sum_{k=1}^p \frac{(\pm
    \lambda_j)^k}{k!} + \sum_{k = p+1}^{d} (\pm \lambda_j)^k
  h_{0,k} \\
  & = & e^{\pm \lambda_j} - t_p \left[ e^{\pm \lambda_j \cdot} \right] (1)
  + \sum_{\ell=1}^r (\pm \lambda_j)^{p+2\ell-1} h_{0,p+2\ell-1} + (\pm
  \lambda_j)^{p+2\ell} h_{0,p+2\ell}.
\end{eqnarray*}
Again, we add and subtract to obtain
\begin{eqnarray*}
  0 & = & \left( e^{\lambda_j} \pm e^{-\lambda_j} \right) - t_p \left[ e^{\lambda_j
      \cdot} \pm e^{-\lambda_j \cdot} \right] (1) \\
  & & + \sum_{\ell=1}^r \left(
    \lambda_j^{p+2\ell-1} \pm (- \lambda_j)^{p+2\ell-1} \right)
  h_{0,p+2\ell-1}  + \left( \lambda_j^{p+2\ell} \pm (- \lambda_j)^{p+2\ell}
  \right) h_{0,p+2\ell}.
\end{eqnarray*}
This again decomposes depending on the parity of $p$. Supposing that
$p$ is even, we get for $j=1,\dots,r$
\begin{eqnarray*}
 \sum_{\ell=0}^{r-1} \lambda_j^{2\ell} \, h_{0,p+2\ell+1} & = & -
  \frac{\left( e^{\lambda_j} - e^{-\lambda_j} \right) - t_p \left[
      e^{\lambda_j \cdot} - e^{-\lambda_j \cdot} \right] (1)}{2
    \lambda_j^{p+1}}, \\
  \sum_{\ell=0}^{r-1} \lambda_j^{2\ell} \, h_{0,p+2\ell+2} & = & -
  \frac{\left( e^{\lambda_j} + e^{-\lambda_j} \right) - t_p \left[
      e^{\lambda_j \cdot} + e^{-\lambda_j \cdot} \right] (1)}{2
    \lambda_j^{p+2}},
\end{eqnarray*}
and since the polynomials $1,x^2,\dots,x^{2r-2}$ form a Chebychev
system on $\RR_+$, this system of equations has a unique
solution. Defining the vectors
\begin{eqnarray*}
w_+ & := & \left[ -\frac{e^{\lambda_j} + e^{-\lambda_j} - t_p \left[
      e^{\lambda_j \cdot} + e^{-\lambda_j \cdot} \right] (1)}{2
    \lambda_j^{p+2}} \;:\; j=1,\dots,r \right], \\
w_- & := & \left[ -\frac{e^{\lambda_j} - e^{-\lambda_j} - t_p \left[
      e^{\lambda_j \cdot} - e^{-\lambda_j \cdot} \right] (1)}{2
    \lambda_j^{p+1}} \;:\; j=1,\dots,r \right],
\end{eqnarray*}
and the \emph{Vandermonde matrices}
$$
L_\Lambda = \left[ \lambda_j^{2\ell} \;:\;
  \begin{array}{c}
    j = 1,\dots,r \\ \ell = 0,\dots,r-1
  \end{array}
\right] \in \RR^{r\times r},
$$
we can therefore write down the construction of the cancellation
operator explicitly.

\begin{lemma}
  The condition (\ref{eq:HLambdaCancel}) can be satisfied by setting
  \begin{eqnarray*}
    \left[ h_{0,p+2\ell+1} \;:\; \ell = 0,\dots,r-1 \right] = \left\{
      \begin{array}{ccl}
        L_\Lambda^{-1} w_-, & \quad & p \in 2 \NN, \\
        L_\Lambda^{-1} w_+, & \quad & p \in 2 \NN+1,
      \end{array}
    \right. \\
    \left[ h_{0,p+2\ell+2} \;:\; \ell = 0,\dots,r-1 \right] = \left\{
      \begin{array}{ccl}
        L_\Lambda^{-1} w_+, & \quad & p \in 2 \NN, \\
        L_\Lambda^{-1} w_-, & \quad & p \in 2 \NN+1.
      \end{array}
    \right.
  \end{eqnarray*}
\end{lemma}

\noindent
The completion of $\cH_{d,\Lambda}$ by means of $\cR$ is now an obvious
extension of (\ref{eq:RzDefinition}), namely
\begin{equation}
  \label{eq:RzGeneralDef}
  R^* (z) = L_{d,\Lambda} \, \Delta_\Lambda^* (z) \, L_{d,\Lambda}^{-1},
\end{equation}
where
\begin{equation}
  \label{eq:RzGeneralDef1}
  L_{d,\Lambda} := \left[
    \begin{array}{ccccc}
      \lambda_1^{p+1} & (-\lambda_1)^{p+1} & \dots & \lambda_r^{p+1} &
      (-\lambda_r)^{p+1} \\
      \vdots & \vdots & \ddots & \vdots & \vdots \\
     \lambda_1^{d} & (-\lambda_1)^{d} & \dots & \lambda_r^{d} &
      (-\lambda_r)^{d}
    \end{array}
  \right] \in \RR^{2r \times 2r}\,,
\end{equation}
and
\begin{eqnarray}
  \nonumber
  \Delta_\Lambda^* (z) & := & \mbox{diag}\, \left[ \Delta_{\pm
      \lambda_j}^* (z) \;:\; j=1,\dots,r \right] \\
  \label{eq:RzGeneralDef2}
  & = & \left[
    \begin{array}{ccccc}
      z^{-1} - e^{\lambda_1} \\
      & z^{-1} - e^{-\lambda_1} \\
      & & \ddots \\
      & & & z^{-1} - e^{\lambda_r} \\
      & & & & z^{-1} - e^{-\lambda_r}
    \end{array}
  \right].
\end{eqnarray}
Since $L_{d,\Lambda}$ is the transpose of a Vandermonde matrix, it is nonsingular.


\section{Factorization}
\label{sec:factorization}
The main result for the use of cancellation operators is related to
the factorization of any subdivision operator that satisfies the
spectral condition.

\begin{theorem}\label{thm:Factorization}
  If the subdivision operator $\cS_{\Ab^{[n]}}$  satisfies the
  $V_{d,\lambda}$-spectral condition, then there exists a mask $\Bb^{[n]}\in \ell_{00}^{(d+1)\times (d+1)}(\ZZ)$
  such that
  \begin{equation}\label{eq:factSA}
    \cH_{d,2^{-(n+1)} \Lambda} \cS_{\Ab^{[n]}} = \cS_{\Bb^{[n]}}
    \cH_{d,2^{-n} \Lambda},
  \end{equation}
  or, in terms of symbols,
  \begin{equation}\label{eq:factSymbolA}
    H_{d,2^{-(n+1)} \Lambda}^*(z) \, \left( A^{[n]} \right)^*(z) = \left(
      B^{[n]} \right)^*(z) \, H_{d,2^{-n} \Lambda}^*(z^2).
  \end{equation}
\end{theorem}

\noindent
In order to prove this theorem, we first give some results about the
factorization of (subdivision and convolution) operators which
annihilate the space $V_{d,\Lambda}$.

\begin{theorem}\label{thm:CFactorization}
  If $\Cb \in \ell^{(d+1) \times (d+1)}_{00} (\ZZ)$ is a finitely supported mask
  such that $\cS_\Cb V_{d,\Lambda} = 0$, then there exists a finitely
  supported mask $\Bb\in \ell^{(d+1)  \times (d+1) }_{00} (\ZZ)$ such that
  $\cS_\Cb = \cS_\Bb \, \cH_{d,\Lambda}$.
\end{theorem}

\begin{pf}
  We first recall from \cite{MerrienSauer12} that whenever $\cS_\cb \Pi_p
  = 0$, then there exists $\Bb \in \ell_{00}^{(d+1)  \times (d+1) } (\ZZ)$
  such that
  $$
  \cS_\Cb = \cS_{\Bb} \left[
    \begin{array}{cc}
      \cT_{p} & 0 \\ 0 & I
    \end{array}
  \right],
  $$
  and $\Bb$ has a symbol with structure
  $$
  B^* (z) = \left[ B_p^* (z), C_{2r}^* (z) \right] := \left[ b_0^*
    (z),\dots,b_p^* (z), c_{p+1}^* (z), \dots, c_d^* (z) \right],
  $$
  where $c^*_{p+1}(z),\dots,c^*_d(z)$ are columns of the original $C^*(z)$. We define the matrix sequence
  $$
  \Wb := \left[ \vb_{e^{\lambda_1 \cdot}}, \vb_{e^{-\lambda_1
        \cdot}},\dots,\vb_{e^{\lambda_r \cdot}}, \vb_{e^{-\lambda_r
        \cdot}} \right] \in \RR^{(d+1) \times 2r}.
  $$
  By assumption, $\cS_\Cb \Wb =0$ and $\cH_{d,\Lambda} \Wb =\left[
      \begin{array}{cc}
        \cT_p & \cQ\\ 0 & \cR
      \end{array}
    \right] \Wb	= 0$, and we thus get
  \begin{eqnarray*}
    0 & = & \cS_\Cb \Wb = \cS_\Bb \left[
      \begin{array}{cc}
        \cT_{p} & 0 \\ 0 & I
      \end{array}
    \right] \Wb = \cS_\Bb \left( \left[
        \begin{array}{cc}
          \cT_{p} & 0 \\ 0 & I
        \end{array}
      \right] - \cH_{d,\Lambda} \right) \Wb \\
    & = & \cS_\Bb \left[
      \begin{array}{cc}
        0 & -\cQ\\ 0 & I - \cR
      \end{array}
    \right] \Wb
    = \cS_\Bb \left[
      \begin{array}{cc}
        0 & -\cQ \\ 0 & I
      \end{array}
    \right] \Wb
    = \cS_\Bb \left[
      \begin{array}{c}
        -\cQ L_{d,\Lambda} \\ L_{d,\Lambda}
      \end{array}
    \right] \mbox{diag}\, \left( e^{\pm \Lambda\cdot} \right) \\
    & = & \sum_{\alpha \in \ZZ} \left( - B_p ( \cdot - 2\alpha ) Q( \cdot - 2\alpha ) +
      C_{2r} (\cdot - 2\alpha) \right) L_{d,\Lambda} \, \mbox{diag} \left( e^{\pm
        \Lambda\cdot} \right),
  \end{eqnarray*}
	where
	$$\mbox{diag} \left( e^{\pm
        \Lambda\cdot} \right):=\left[
				\begin{array}{ccccc}
				e^{\lambda_1\cdot} \\
				& 	e^{-\lambda_1\cdot}  \\
				&& \ddots \\
				&&& 	e^{\lambda_r\cdot} \\
					&&&& 	e^{-\lambda_r\cdot} \end{array}
					\right].
					$$
  This implies that for $\epsilon \in \{0,1\}$ and $j=1,\dots,r$ we
  must have
  \begin{eqnarray}
    \label{eq:thmFactorizationPf2a}
    0 & = & e^{\lambda_j \cdot} \sum_{\alpha \in \ZZ} \left( -B_p
      (\epsilon+2\alpha) Q (\epsilon + 2\alpha ) + C_{2r}
      (\epsilon+2\alpha) \right) L_{d,\Lambda} e_{2j-1} e^{-\lambda_j
      \alpha}, \\
    \label{eq:thmFactorizationPf2b}
    0 & = & e^{-\lambda_j \cdot} \sum_{\alpha \in \ZZ} \left( -B_p
      (\epsilon+2\alpha) Q( \epsilon + 2\alpha ) + C_{2r}
      (\epsilon+2\alpha) \right) L_{d,\Lambda} e_{2j} e^{\lambda_j \alpha},
  \end{eqnarray}
  with $e_j$ the standard $j$-th unit vector in $\RR^{d+1}$, from which it follows that
  $$
  0 = \left( -B_p^* Q^* + C^*_{2r} \right) L_{d,\Lambda} e_{2j-1} \left( \pm
    e^{-\lambda_j/2} \right)
  = \left( -B_p^* Q^*+ C^*_{2r} \right) L_{d,\Lambda} e_{2j} \left( \pm
    e^{\lambda_j/2} \right),
  $$
  hence,
  \begin{equation}
    \label{eq:thmFactorizationPf1}
    \left( -B_p^* (z) Q^*(z^2) + C^*_{2r} (z) \right) L_{d,\Lambda} e_{2j-1}
    = \left( z^{-2} - e^{\lambda_j} \right) \, b^*_{2j-1} (z), \qquad
    j=1,\dots,r,
  \end{equation}
  and
  \begin{equation}
    \label{eq:thmFactorizationPf2}
    \left( -B_p^* (z) Q^*(z^2)+ C^*_{2r} (z) \right) L_{d,\Lambda} e_{2j} =
    \left( z^{-2} - e^{-\lambda_j} \right) \, b^*_{2j} (z), \qquad
    j=1,\dots,r.
  \end{equation}
  Setting $B_{2r}^* (z)= \left[ b_j^*(z) \;:\; j=1,\dots,2r \right]$,
  (\ref{eq:thmFactorizationPf1}) and (\ref{eq:thmFactorizationPf2})
  can be conveniently combined into
  $$
  \left( -B_p^* (z) Q^*(z^2) + C_{2r}^* (z) \right) L_{d,\Lambda} = B_{2r}^* (z) \,
  \Delta_{\Lambda}^* (z^2)
  $$
  which leads to
  $$
  C_{2r}^* (z) = B_{2r}^* (z) \, L_{d,\Lambda}^{-1} L_{d,\Lambda}
  \Delta_\Lambda^* (z^2) L_{d,\Lambda}^{-1} + B_p^* (z) Q^*(z^2),
  $$
  and consequently
  \begin{eqnarray*}
    B^* (z) & = & \left[ B_p^* (z), C_{2r}^* (z) \right]
    = \left[ B_p^* (z), B_{2r}^* (z) L_{d,\Lambda}^{-1} \,
      L_{d,\Lambda} \Delta_{\Lambda}^* (z^2) L_{d,\Lambda}^{-1} + B_p^*
      (z) Q^*(z^2) \right] \\
    & = & \left[ B_p^* (z), B_{2r}^* (z) L_{d,\Lambda}^{-1} \right] \left[
      \begin{array}{cc}
        I & Q^*(z^2) \\ 0 & R^* (z^2)
      \end{array}
    \right].
  \end{eqnarray*}
  This eventually gives
  \begin{eqnarray*}
    C^* (z) & = & B^* (z) \left[
      \begin{array}{cc}
        T_{p}^* (z^2) & 0 \\ 0 & I
      \end{array}
    \right] \\
    & = & \left[ B_p^* (z), B_{2r}^* (z) L_{d,\Lambda}^{-1} \right] \left[
      \begin{array}{cc}
        I & Q^*(z^2) \\ 0 & R^* (z^2)
      \end{array}
    \right] \left[
      \begin{array}{cc}
        T_{p}^* (z^2) & 0 \\ 0 & I
      \end{array}
    \right] \\
    & = & \left[ B_p^* (z), B_{2r}^* (z) L_{d,\Lambda}^{-1} \right] \left[
      \begin{array}{cc}
        T_{d-2}^* (z^2) & Q^*(z^2)\\
        0 & R^* (z^2)
      \end{array}
    \right] \\
    & = & \left[ B_p^* (z), B_{2r}^* (z) L_{d,\Lambda}^{-1} \right]
    H_{d,\Lambda}^* (z^2)\,,
  \end{eqnarray*}
  and completes the proof.
\end{pf}

\noindent
As a consequence of Theorem~\ref{thm:CFactorization} and
Remark~\ref{rem_marianto} we get the desired result that extends the
observations made in the introduction.

\begin{corollary}\label{cor:CFactorization}
  If $\Cb^{[n]} \in \ell^{(d+1) \times (d+1)}_{00} (\ZZ)$ is 
  such that $\cS_{\Cb^{[n]}} \vb_{f,n}= 0$, $f\in V_{d,\Lambda}$, then
  there exists a finitely supported mask $\Bb^{[n]}\in \ell^{(d+1)
    \times (d+1)}_{00} (\ZZ)$ such that $\cS_{\Cb^{[n]}} =
  \cS_{\Bb^{[n]}} \cH_{d,2^{-n}\Lambda}$.
\end{corollary}

\noindent
Using this result, Theorem~\ref{thm:Factorization} is now easy to prove.
\bigskip

\begin{pfof}{Theorem~\ref{thm:Factorization}}
  Set $\cS_{\Cb^{[n]}} := \cH_{d,2^{-n-1}\Lambda} \,
  \cS_{\Ab^{[n]}}$. Since for $f\in V_{d,\Lambda}$ we have
  $$
  \cS_{\Cb^{[n]}} \vb_{f,n} = \cH_{d,2^{-n-1}\Lambda} \,
  \cS_{\Ab^{[n]}}\vb_{f,n} = \cH_{d,2^{-n-1}\Lambda} \, \vb_{f,{n+1}}
  = 0,
  $$
  it follows from Corollary~\ref{cor:CFactorization} that there exists
  $\Bb^{[n]}$ such that
  $$
  \cH_{d,2^{-n-1}\Lambda} \, \cS_{\Ab^{[n]}} = \cS_{\Bb^{[n]}} \,
  \cH_{d,2^{-n}\Lambda},
  $$
  as claimed.
\end{pfof}

\begin{remark}
  Note that the factorization (\ref{eq:factSA}) of $\cS_{\Ab^{[n]}}$
  is equivalent to the following factorization of $\cS_{\widetilde
    \Ab^{[n]}}$:
 \begin{equation}\label{eq:factStildeA}
   \widetilde \cH_{d,\Lambda}^{[n+1]} \cS_{\widetilde \Ab^{[n]}} =
   \cS_{\widetilde \Bb^{[n]}} \widetilde \cH_{d,\Lambda}^{[n]},
 \end{equation}
 where
 $$
 \widetilde \Hb_{d,\Lambda}^{[n]} := D^n \widetilde
 \Hb_{d,2^{-n}\Lambda} D^{-n},\qquad \widetilde \Bb^{[n]} :=D^{-n-1}
 \widetilde \Bb^{[n]} D^{n}.
 $$
\end{remark}

\noindent
A careful inspection of the  proof of Theorem~\ref{thm:CFactorization}
shows that the factorization can also be extended to convolution
operators.

\begin{theorem}\label{thm:ConvFactorization}
  If $\Cb\in  \ell^{(d+1) \times (d+1)}_{00} (\ZZ)$ is 
  such that $\Cb * V_{d,\Lambda} =  0$, then there exists a finitely
  supported mask $\Bb\in \ell^{(d+1) \times (d+1)}_{00} (\ZZ)$ 
	such that $\Cb = \Bb * \Hb_{d,\Lambda}$.
\end{theorem}

\begin{pf}
  The proof follows exactly the lines of the one of
  Theorem~\ref{thm:Factorization} except that
  (\ref{eq:thmFactorizationPf2a}) and (\ref{eq:thmFactorizationPf2b})
  become
  \begin{eqnarray*}
    0 & = & e^{\lambda_j \cdot} \sum_{\alpha \in \ZZ} \left( -B_p (\alpha)Q (\alpha) +
      C_{2r} (\alpha) \right) L_{d,\Lambda} e_{2j-1} e^{-\lambda_j
      \alpha}, \qquad j=1,\dots,r, \\ 
    0 & = & e^{-\lambda_j \cdot} \sum_{\alpha \in \ZZ} \left( -B_p (\alpha) Q(\alpha) +
      C_{2r} (\alpha) \right) L_{d,\Lambda} e_{2j} e^{\lambda_j
      \alpha}, \qquad j=1,\dots,r, 
  \end{eqnarray*}
  that is,
  \begin{eqnarray*}
    \left( -B_p^* (z) Q^*(z) + C^*_{2r} (z) \right) L_{d,\Lambda} e_{2j-1}
    & = & \left( z^{-1} - e^{\lambda_j} \right) \, b^*_{2j-1} (z), \qquad
    j=1,\dots,r, \\
    \left( -B_p^* (z) Q^*(z) + C^*_{2r} (z) \right) L_{d,\Lambda} e_{2j} &
    = & \left( z^{-1} - e^{-\lambda_j} \right) \, b^*_{2j} (z), \qquad
    j=1,\dots,r.
  \end{eqnarray*}
  From there on the arguments can be repeated literally to yield that
  \begin{equation}
    \label{eq:ConvOpFactor}
    C^* (z) = B^* (z) \, H_{d,\Lambda}^* (z).
  \end{equation}
  Finally, observe that in the same way the argument from
  \cite{MerrienSauer12} can be modified to give the initial
  factorization by means of the Taylor operator.
\end{pf}

\medskip \noindent
Since $\cH_{d,\Lambda}$ is a convolution operator itself and since
(\ref{eq:ConvOpFactor}) can be reformulated as the fact that for
\emph{any} $\Cb$ that annihilates $V_{d,\Lambda}$, the
Laurent polynomial $\det C^* (z)$ must be divisible by $\det
H_{d,\Lambda}^* (z)$, this operator is a particular annihilator of
$V_{d,\Lambda}$.

\begin{corollary}\label{cor:SpecialAnni}
  The operator $\cH_{d,\Lambda}$ is a minimal annihilator for
  $V_{d,\Lambda}$.
\end{corollary}

\begin{corollary}
  The Taylor operator $\cT_d$ is a minimal annihilator for
  $V_{d,\emptyset}$.
\end{corollary}

\section{Examples}
To illustrate the results of the preceding sections, we construct
two matrix subdivision schemes which reproduce, by construction,
polynomials and exponential from the spaces
$$
V_{2,\lambda} = \mbox{\rm span} \left\{ 1,e^{-\lambda x},e^{\lambda x}
\right\}, \qquad V_{3,\lambda} = \mbox{\rm span}
\left\{1,x,e^{-\lambda x},e^{\lambda x} \right\},
$$
and explicitly verify for these cases the factorization property via
the annihilators in (\ref{esempio:H_2}) and (\ref{esempio:H_3}).

To construct the first vector Hermite subdivision scheme, we start
with a sufficiently smooth real valued function $f$ and define the
initial sequence of vector data $\pb^0 =\left( p (\alpha) := [f(\alpha),
  f'(\alpha), f'' (\alpha)]^T \;:\; \alpha \in \ZZ \right)$ from which
we construct in each interval the functions $g^0_\alpha,
g^0_{\alpha+1} \in V_{2,\Lambda}$ such that they solve Hermite
interpolation problems at $\alpha$ and $\alpha+1$ based on the data
$p^0 (\alpha)$ and $p^0 (\alpha+1)$, respectively. It is easy to
verify that these interpolation problems admit a unique solution in
$V_{2,\lambda}$. This leads to the general interpolatory subdivision rules
\begin{equation}
  \label{def:sub_rulesESEMPIO1}
  \begin{array}{rcl}
    p^{n+1} (2\alpha) & = & p^{n} (\alpha),\\\noalign{\medskip}
    p^{n+1} (2\alpha+1) & = & \frac 12 \left( g^n_\alpha \left(
        \frac{2\alpha+1}{2^{n+1}} \right) + g^n_{\alpha+1} \left(
        \frac{2\alpha+1}{2^{n+1}} \right) \right),
  \end{array}
  \qquad n \in \NN_0.
\end{equation}
It turns out that matrix masks of the interpolatory Hermite
subdivision scheme defined as in (\ref{eq:HermSubd0}) consist of three
nonzero $3\times 3$ matrices. 
The symbol of the scheme at the $n$-th iteration is
{\footnotesize
  \begin{equation}\label{def:Aes1} (A^{[n]})^*(z)=
    \frac 1{16z}\left[
      \begin {array}{ccc}
        8\, \left( z+1 \right) ^{2}&
        \displaystyle{\frac {4}{\lambda_n} \left( z^2-1\right)
          \sinh \frac{\lambda_n}{2} }&\displaystyle{\frac {4}{\lambda_n^{2}}
          \left( 1+{z}^{2} \right)  \left( \cosh
            \frac{\lambda_n}{2}-2 \right)} \\\noalign{\medskip} 
        0&2(1+z^2)\cosh  \frac{\lambda_n}{2}+8\,z&
        \displaystyle{\frac {2}{\lambda_n}
          \left( z^2-1 \right)  \ \sinh  \frac{\lambda_n}{2} }
        \\\noalign{\medskip}0&\lambda_n\, \left( z^2 -1\right)
        \sinh  \frac{\lambda_n}{2}
        &
        (1+z^2)\ \cosh  \frac{\lambda_n}{2}+4\,z
      \end {array} \right],
  \end{equation}
}
with the abbreviation $\lambda_n:={2}^{-n}\lambda$.

\begin{figure}
\begin{center}
\includegraphics[width=10cm]{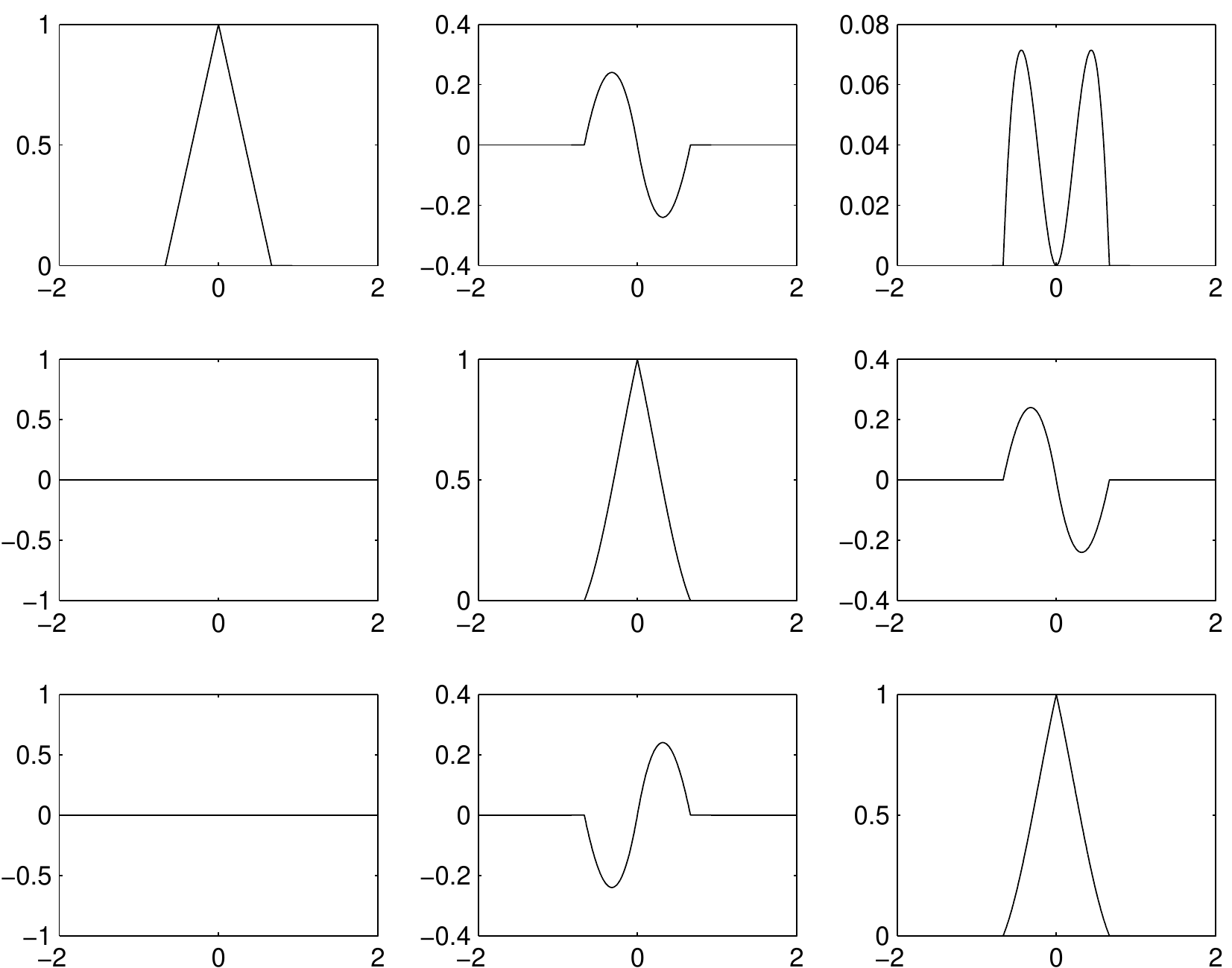}
	\caption{Result after $12$ iterations of the $3\times 3$ non-stationary subdivision scheme in  (\ref{def:sub_rulesESEMPIO1}).}
\label{fig:A3x3}
\end{center}
\end{figure}

Observe that the determinant of $(A^{[n]})^*(z)$, $n \in \NN_0$, factorizes into
$$
\det (A^{[n]})^*(z) = \frac{(z+1)^2 e^{-\lambda_n}
  \left(e^{\frac{\lambda_n}{2}}+z\right)^2 \left(z
    e^{\frac{\lambda_n}{2}}+1\right)^2}{64 z^3}.
$$
The resulting subdivision scheme appears to be convergent since, when
starting the subdivision iterations by applying column-wise the
subdivision rules to the delta matrix sequence  the result after
$12$ iterations stabilizes on the matrix function shown in
Fig.~\ref{fig:A3x3}, but a specific convergence analysis is not in the
scope of this paper. 

\begin{figure}
\begin{center}
\includegraphics[width=10cm]{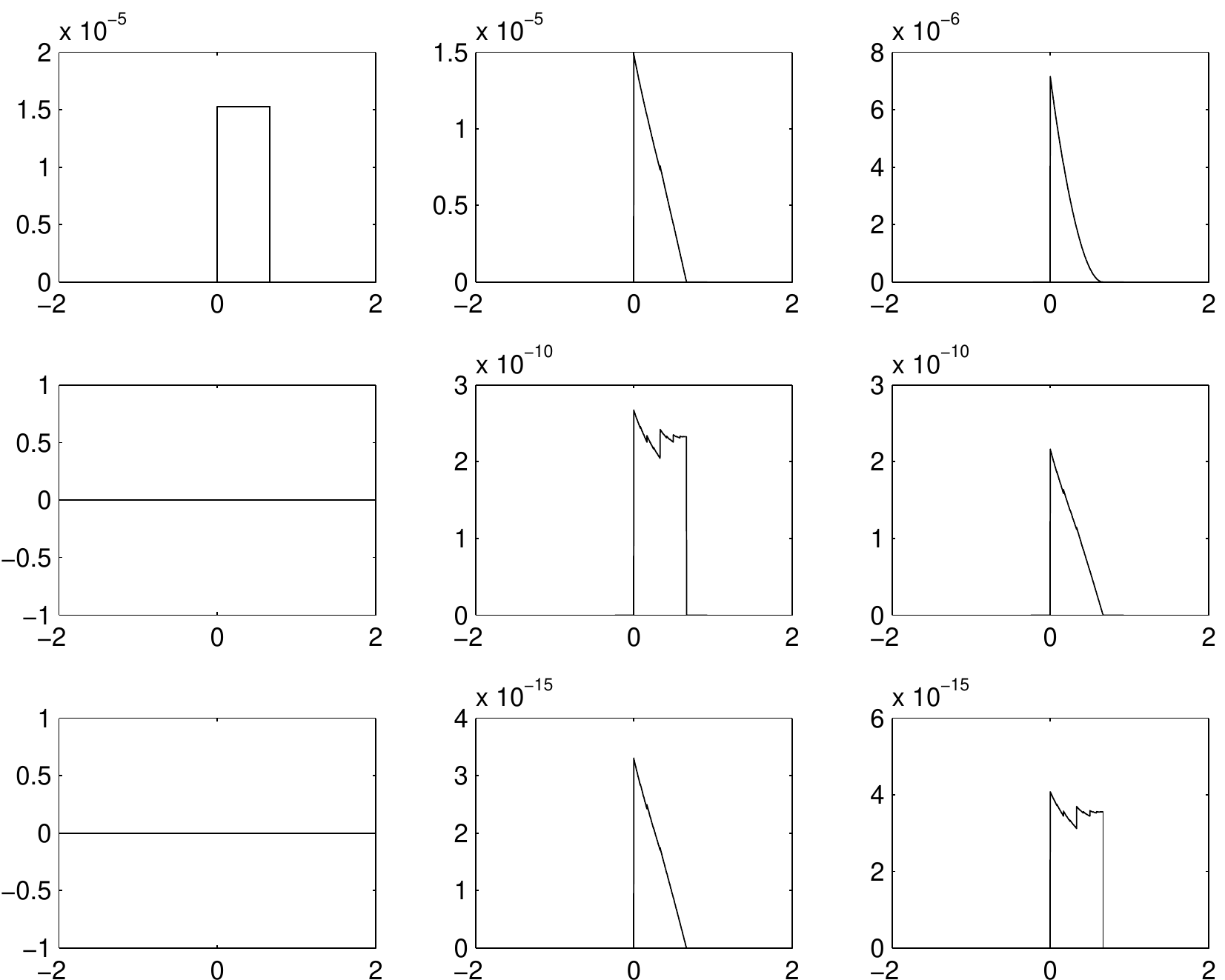}
	\caption{Result after $12$ iterations of the $3\times 3$ non-stationary subdivision scheme based on (\ref{def:Bes1}).}\label{fig:B3x3}
\end{center}
\end{figure}

By construction this scheme satisfies the $V_{2,\lambda}$-spectral
condition and according to Theorem~\ref{thm:Factorization} it is
possible to find a subdivision operator $\cS_{\Bb^{[n]}}$ such that
the factorization (\ref{eq:factSA}) holds true. At the $n$-th
iteration, its symbol is given by: 
\begin{equation}\label{def:Bes1}
  \left( B^{[n]} \right)^*(z)=\frac{1}{16}
  \left[
    \begin {array}{ccc} 8+8\,z&-4\,\displaystyle{\frac {\sinh
          \frac{\lambda_n}{2}}{\lambda_n}}&4\,\displaystyle{\frac
        {\cosh  \frac{\lambda_n}{2}-2}{{\lambda_n}^{2}}} 
      \\\noalign{\medskip}0&2\,\displaystyle{\cosh
        \frac{\lambda_n}{2}+4\,z}&-2\,\displaystyle{\frac {\sinh
          \frac{\lambda_n}{2}}{\lambda_n}}\\\noalign{\medskip}0&-\lambda_n\, 
      \sinh  \frac{\lambda_n}{2}&
      \cosh  \frac{\lambda_n}{2}+2\,z\end {array} \right].
\end{equation}
\noindent The corresponding subdivision scheme seems to be
zero-convergent, see Figure~\ref{fig:B3x3}, and hence contractive as
one should expect.

To construct the second example, we define the initial sequence of
vector data $\pb^0= \left( [f(\alpha), f'(\alpha), f''(\alpha),
  f'''(\alpha)]^T \;:\; \alpha \in \ZZ \right)$ and apply the same
construction as above, just in $V_{3,\Lambda}$.

\begin{figure}
\begin{center}
\includegraphics[width=10cm]{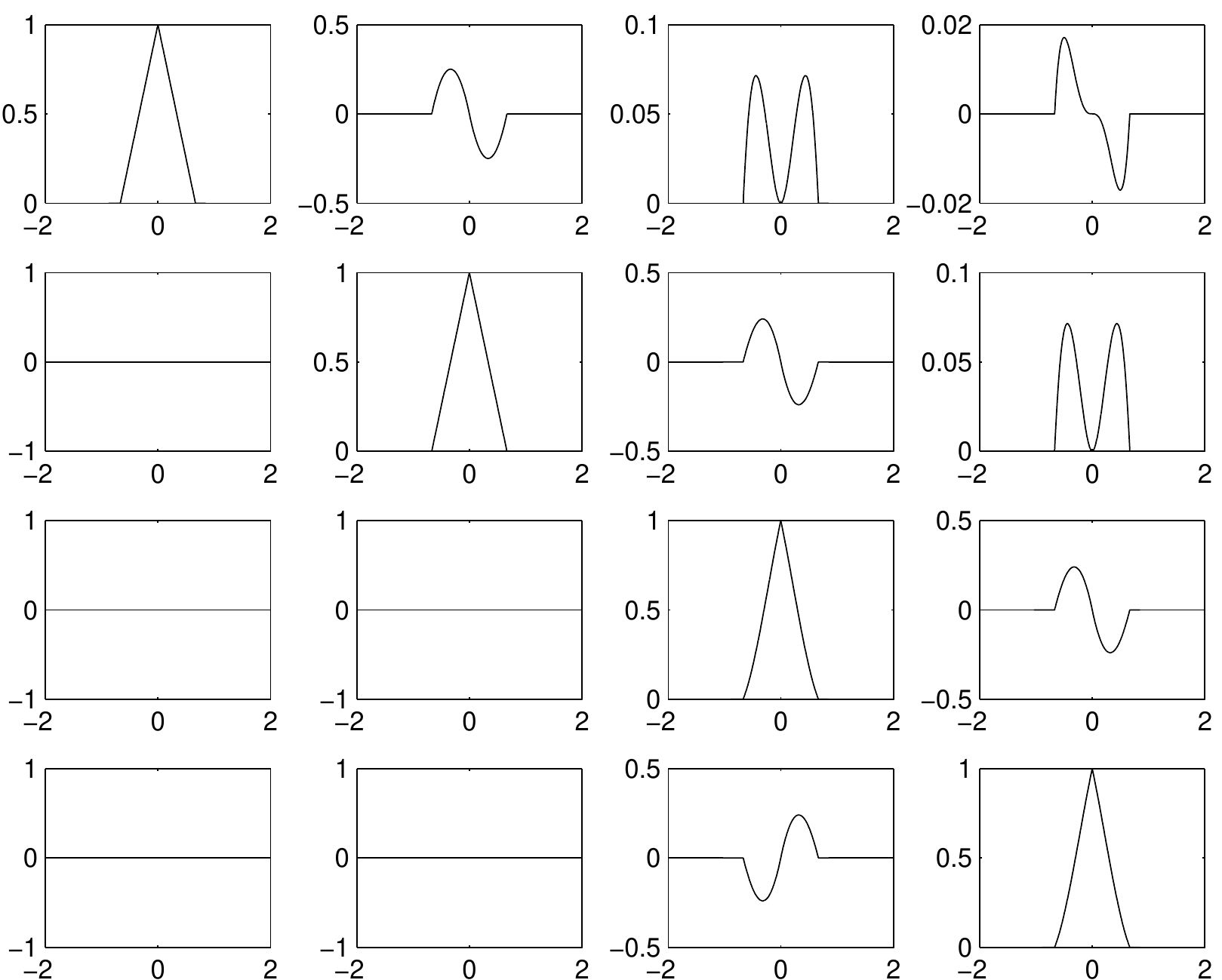}
	\caption{Result after $12$ iterations of the $4\times 4$ non-stationary subdivision scheme based on (\ref{def:Aes2}).}\label{fig:A4x4}
\end{center}

\end{figure}

The symbol at level $n$ can be computed explicitly as
\begin{equation}\label{def:Aes2}
  {\scriptsize  
    \left( A^{[n]} \right)^*(z)=\frac 1{32 z}
    \left[ \begin {array}{cccc} 16\, \left( 1+z \right) ^{2}&
        8\, \left(z^2-1 \right) &\displaystyle{\frac {8}{\lambda_n^2}\, \left( 1+{z}^{2} \right)
          \left(\cosh  \frac{\lambda_n}{2}-2 \right)
        }&
        \displaystyle{\frac {8}{\lambda^3}\, { \left( z^2-1\right)
            \left( \lambda+\sinh  \frac{\lambda_n}{2}
            \right) }}\\\noalign{\medskip}
        0&8\, \left( z+1 \right) ^{2}&
        \displaystyle{\frac {4}{\lambda_n} \left( z^2-1\right)   \sinh  \frac{\lambda_n}{2} }
        &\displaystyle{\frac {4}{\lambda_n^{2}}
          \left( 1+{z}^{2} \right)  \left(\cosh  \frac{\lambda_n}{2}-2 \right)} \\\noalign{\medskip}
        0& 0&2(1+z^2)\ \cosh  \frac{\lambda_n}{2}+8\,z&
        \displaystyle{\frac {2}{\lambda_n}
          \left( z^2-1 \right)  \ \sinh  \frac{\lambda_n}{2} }
        \\\noalign{\medskip}
        0& 0&\lambda_n\, \left( z^2-1 \right)
        \sinh  \frac{\lambda_n}{2} &
        (1+z^2)\cosh  \frac{\lambda_n}{2}+4\,z
      \end {array} \right],}
\end{equation}
The determinant of $\left( A^{[n]} \right)^*(z)$ factorizes into
$$
\det \left( A^{[n]} \right)^*(z) = \frac{(z+1)^4 e^{-\lambda_n}
  \left(e^{\frac{\lambda_n}{2}}+z\right)^2 \left(z
    e^{\frac{\lambda_n}{2}}+1\right)^2}{1024 \, z^4}.
$$
Evidence for the convergence of this scheme is given in
Fig.~\ref{fig:A4x4}, where we show the plot of $12$ iterations of the
scheme applied to the delta matrix sequence.

\begin{figure}
\begin{center}
\includegraphics[width=10cm]{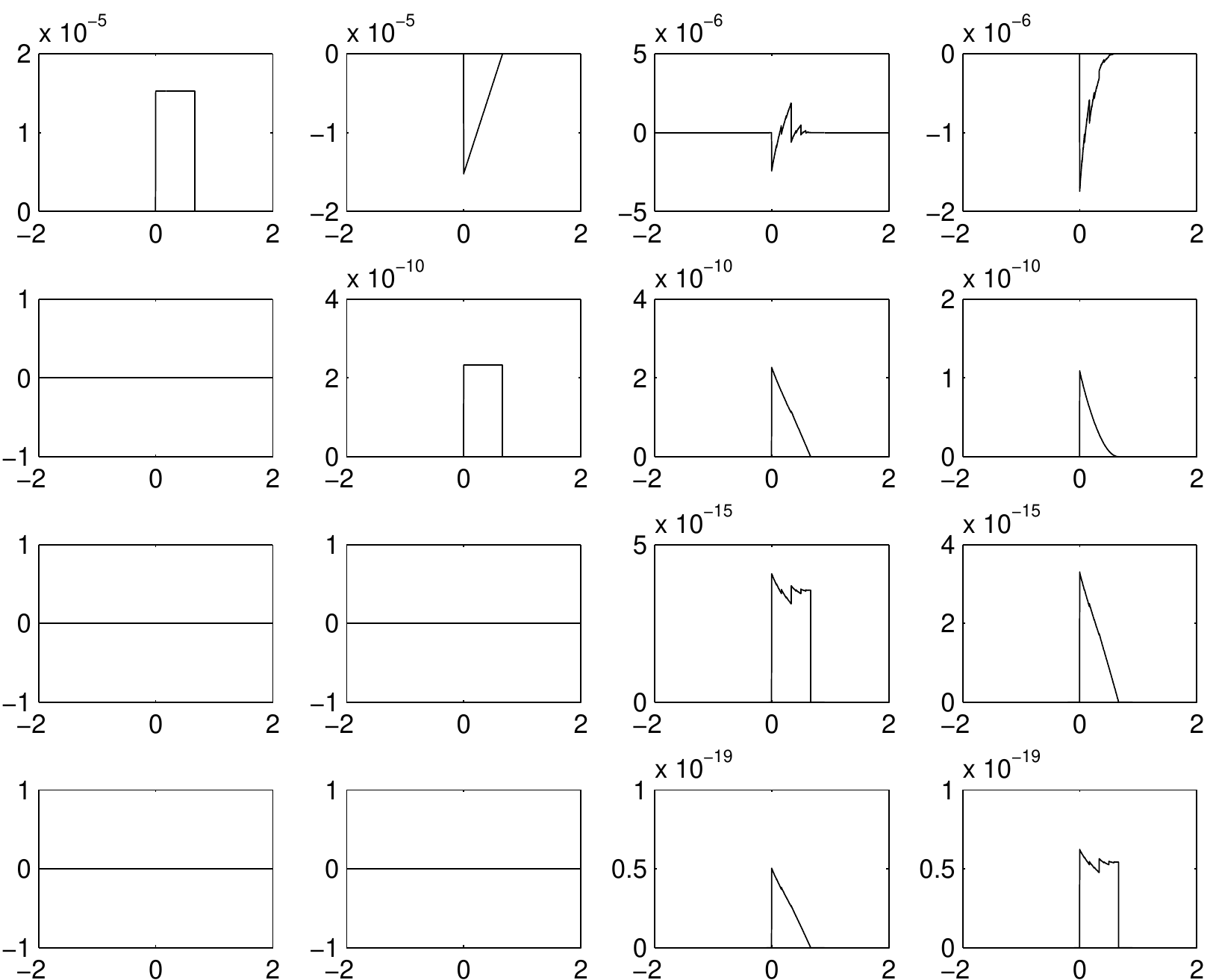}
	\caption{Result after $12$ iterations of the $3\times 3$ non-stationary subdivision scheme based on (\ref{def:Bes2}).}\label{fig:B4x4}
\end{center}
\end{figure}

This scheme satisfies the $V_{3,\lambda}$-spectral condition  and
therefore admits the  factorization  (\ref{eq:factSA}) with
\begin{equation}\label{def:Bes2}
  \left( B^{[n]} \right)^*(z)
  =\frac 1{32}\left[ \begin {array}{cccc}
      16+16\,z&-8&8\,\displaystyle{\frac {\cosh
          \frac{\lambda_n}{2}-2}{{\lambda_n}^{2}}}&8\,\displaystyle{\frac
        { 
          \lambda-\sinh  \frac{\lambda_n}{2}}{{\lambda_n}^
          {3}}}\\\noalign{\medskip}0&8+8\,z&-4\,\displaystyle{\frac
        {\sinh
          \frac{\lambda_n}{2}}{\lambda_n}}&4\,\displaystyle{\frac
        {\cosh  \frac{\lambda_n}{2}-2}{{\lambda_n}^{2}}} 
      \\\noalign{\medskip}0&0&2\,\cosh
      \frac{\lambda_n}{2}+4\,z&-2\,\displaystyle{\frac {\sinh
          \frac{\lambda_n}{2}}{\lambda_n}}\\ \noalign{\medskip} 0 & 0
      & \displaystyle{\lambda_n\,
        \sinh  \frac{\lambda_n}{2}}&\displaystyle{\cosh
        \frac{\lambda_n}{2}+2\,z}\end {array} \right],
\end{equation}
which again seems to be a contraction, see Figure~\ref{fig:B4x4}.

We conclude this section by observing that, as $n$ goes to infinity,
the symbols $\left( A^{n]} \right)^*(z)$ in (\ref{def:Aes1}) and
  (\ref{def:Aes2}) tend to
$$
\left( A^{[\infty]} \right)^*(z) = \frac 1{16z}\left[ \begin {array}{ccc}
    8\, \left( z+1 \right) ^{2}
    &4\, \left( z^2-1 \right) &({z}^{2}+1)
    \\\noalign{\medskip}
    0&4\, \left( z+1 \right) ^{2}&2\, \left( z^2-1
    \right)  \\\noalign{\medskip}0&0&2\, \left( z+1
    \right) ^{2}\end {array} \right],
$$
and
$$
\left( A^{[\infty]} \right)^*(z) = \frac 1{96z}\left[ \begin
    {array}{cccc} 48\, \left( z+1 \right) ^{2}&24\, \left( z^2-1
    \right) &6\,({z}^{2}+1)& \left( z^2-1 \right) 
    \\\noalign{\medskip}0&24\, \left( z+1 \right) ^{2}
    &12\, \left( z^2-1 \right) &3\,({z}^{2}+1)
    \\\noalign{\medskip}0&0&12\, \left( z+1 \right) ^{2}&6\, \left( z^2-1
    \right)  \\\noalign{\medskip}0&0&0&6\, \left( z+1 
    \right) ^{2}\end {array} \right],
$$
respectively. These are the symbols of Hermite schemes satisfying a
polynomial space spectral condition. In particular, they  reproduce
polynomials up to the degree  2 and 3, respectively.


\bibliographystyle{amsplain}

\begin{thebibliography}{10}


\bibitem{CharinaContiRomani13} {M.~Charina, C.~Conti, L.~Romani},
  \emph{Reproduction of exponential polynomials by multivariate
    non-stationary subdivision schemes with a general dilation
    matrix}, Numer. Math. (2013), in print.
DOI 10.1007/s00211-013-0587-8.


\bibitem{ContiRomani11}{C.~Conti, L.~Romani}, \emph{Algebraic
    conditions on non-stationary subdivision symbols for exponential
    polynomial reproduction}, J. Comput. Appl. Math. \textbf{236},
  (2011), 543--556.


\bibitem{ContiMerrienRomani14}{C.~Conti, J.L.~Merrien, L.~Romani},
  \emph{Dual Hermite Subdivision Schemes of de Rham-type}, preprint.

\bibitem{DubucMerrien06}
S.~Dubuc and J.-L. Merrien, \emph{Convergent vector and {H}ermite
  subdivision schemes}, Constr. Approx. \textbf{23} (2006), 1--22.

\bibitem{DubucMerrien09}
S.~Dubuc and J.-L. Merrien, \emph{Hermite subdivision schemes and
  {T}aylor polynomials}, Constr. Approx. \textbf{29} (2009),
219--245.

\bibitem{Hamming98}
R.~W. Hamming, \emph{Digital filters}, Prentice--Hall, 1989, Republished by
  Dover Publications, 1998.
	
\bibitem{MerrienSauer12}
J.-L. Merrien and T.~Sauer, \emph{From {H}ermite to stationary subdivision
  schemes in one and several variables}, Advances Comput. Math. \textbf{36}
  (2012), 547--579.
\bibitem{Micchelli96}
C.~A. Micchelli, \emph{Interpolatory subdivision schemes and wavelets}, J.
  Approx. Theory \textbf{86} (1996), 41--71.
\bibitem{UnserBlu05}
  M.~Unser and Th.~Blu, \emph{Cardinal Exponential Splines: Part I ---
    Theory and Filtering Algorithms}, IEEE
  Trans. Sig. Proc. \textbf{53} (2005), 1425--1438.
\end{thebibliography}

\end{document}